\documentclass{article}

\usepackage{amsmath}
\usepackage{graphicx}
\usepackage{amsfonts}
\usepackage{amsthm}
\usepackage{bbm}
\usepackage[colorlinks=true, allcolors=blue]{hyperref}

\newtheorem{theorem}{Theorem}
\newtheorem{lemma}{Lemma}
\newtheorem{corollary}{Corollary}
\usepackage[parfill]{parskip}

\title{Omega theorem for fractional sigma function}
\author{Yuan Qiu, Alexander B. Kalmynin}
\date{}

\begin{document}
\maketitle

\begin{abstract}
The research in the subfield of analytic number theory around error term of summation of sigma functions possesses a history which can be dated back to the mid-19th century when Dirichlet provided an $O(\sqrt{n})$ estimation of error term of summation of $d(n)$ \cite{montgomery2007multiplicative}. Later, G. Voronoi \cite{ivic2012riemann}, G. Kolesnik \cite{kolesnik1982order}, and M.N. Huxley \cite{huxley2003exponential} (to name just a few) contributed more on the upper bound on the error term of summation of sigma functions. As for $\Omega$-theorems, G.H. Hardy was the first contributor. Later researchers on this topic include G.H. Hardy \cite{hardy1929introduction} and T.H. Gronwall\cite{gronwall1913some}, but the amount of academic effort is much sparser than $O$-theorems. This research aims to provide a better $\Omega$-bound for the error term of summation of fractional sigma function $\sigma_{\alpha}(n)$ on the range $0 < \alpha < \frac{1}{2}$, obtaining the result $\Omega((x\ln x)^{\frac{1}{4}+\frac{\alpha}{2}})$. 
\end{abstract}

\section{Introduction}

Let $\sigma_{\alpha}(n)$ denote the fractional sigma function: $\sigma_{\alpha}(n) = \sum\limits_{d|n} d^{\alpha}$. When $0 < \alpha < 1$, it can be shown by summation techniques that  $\sum\limits_{n\leq x} \sigma_{\alpha}(n) = \frac{\zeta(\alpha+1)}{\alpha+1} x^{\alpha+1} + O(x)$ \cite{apostol2013introduction}. Nevertheless, fractional sigma function has not been a hotbed for research in bounding theorems. In fact, the $O(x)$ estimation has been the best known result.

Most research in this field has been diverted to improving the bounding theorem of $d(n)$, which is also called "Dirichlet divisor problem". The problem can be formulated as the estimation of $\Delta(x)$ in the well-known asymptotic formula $\sum_{n\leq x} d(n) = x\log x + (2\gamma-1)x + \Delta(x)$.

In 1849, J.P.G. Dirichlet proved that $\Delta(x) = O(\sqrt{x})$ with the Dirichlet hyperbola method \cite{montgomery2007multiplicative}. In the following years, mathematicians made efforts to improve the error term $O(\sqrt{x})$. In 1904, G. Vonoroi improved the bound to $\Delta(x) = O(x^{1/3}\log x)$ \cite{ivic2012riemann}. In 1982, G. Kolesnik made an improvement to $x^{35/108+\varepsilon}$ \cite{kolesnik1982order}. The best known bound is $\Delta(x) = O(x^{131/416 +\varepsilon})$ due to M.N. Huxley \cite{huxley2003exponential}. Omega-theorem of $\Delta(x)$ is also a direction of research. In 1915, G.H. Hardy proved that $\Delta(x) = \Omega(x^{\frac{1}{4}})$ \cite{hardy1917dirichlet}. There is a widely propagated conjecture that $\Delta(x) = O(x^{1/4 + \varepsilon})$.

Bounding theorem for summation of $\sigma(n)=\sigma_1(n)$ have also been another target of research. From elementary summation techniques a basic O-theorem of summation of $\sigma(n)$ can be obtained: $\sum_{n\leq x} \sigma(n) = \frac{\pi^2 x^2}{12} + O(x\log x)$ \cite{hardy1929introduction}. The best known $O$-bound is obtained by A. Walfisz: $\sum\limits_{n\leq x} \sigma(n) = \frac{\pi^2 x^2}{12} + O\left(x(\log x)^{\frac{2}{3}}\right)$. As for $\Omega$-theorems, In 1913, T.H. Gronwall proved the omega theorem $\Omega(x \log \log x)$ for $\sigma(n)$ \cite{gronwall1913some}. In 1987, P{\'e}termann improved the bound to $\Omega_{-}(x \log \log x)$ \cite{petermann1987omega}. 

Our research primarily aims to produce an $\Omega$-theorem for $\sum\limits_{n\leq x} \sigma_{\alpha}(n)$ when $0 < \alpha < \frac{1}{2}$.

\section{Main results}

We provide a $\Omega$-bound for the remainder term of summation of $\sigma_{\alpha}(n)$, which is better than any previously known bounds.

Our main goal is to prove the following theorem:

\begin{theorem} \label{thm:main}
    Let $\sigma_{\alpha}(n)$ denote the fractional sigma function: $\sigma_{\alpha}(n) = \sum\limits_{d|n} d^{\alpha}$, $0 < \alpha < \frac{1}{2}$. Then
    \begin{equation}
        \sum_{n \leq x} \sigma_{\alpha}(n) = \frac{x^{1+\alpha}}{1+\alpha}\zeta(1+\alpha)+x\zeta(1-\alpha) + E_{\alpha}(x)
    \end{equation}
    where $E_\alpha(x) = \Omega\left((x\ln x)^{\frac{1}{4}+\frac{\alpha}{2}}\right)$.
\end{theorem}

Before proving this result, we will show that $E_{\alpha}$ can be approximated up to $O(x^{\frac{1}{4}+\frac{\alpha}{2}})$ by a certain sum of oscillating terms.

\begin{theorem} \label{E_alpha}
    For any $x^{\frac{1}{2}+\alpha+\varepsilon} \ll N \ll x^{\frac{1}{2\alpha}-\varepsilon}$, 
\begin{equation}
    E_{\alpha}(x) = \frac{x^{\frac{1}{4}+\frac{\alpha}{2}}}{\sqrt{2}\pi}\sum_{m \leq N} \sigma_{\alpha}(m)m^{-\frac{3}{4}-\frac{\alpha}{2}}\cos\left(4\pi\sqrt{mx}-\frac{\pi}{4}\right) + O(x^{\frac{1}{4}+\frac{\alpha}{2}})
\end{equation}
    
\end{theorem}

Note that in the above theorem and throughout the paper, $\varepsilon$ indicates a arbitrarily small positive number.

In the process of proving Theorem \ref{E_alpha}, we also reveal the following bound of $E_{\alpha}(x)$:

\begin{corollary} \label{intermediate_bound}
    \begin{equation}
        E_{\alpha}(x) = O\left(x^{\frac{2\alpha^2+3\alpha+1}{2\alpha+3}+\varepsilon}\right)
    \end{equation}
\end{corollary}

This bound is better than $O(x)$ when $\alpha < \frac{\sqrt{17}-1}{4}$.

Note that we are aware of the following "combinatorial" $\Omega$-bound of $E_{\alpha}(x)$, namely

\begin{theorem}
    For $x\to +\infty$ we have
    \begin{equation}
        \sup_{n\leq x} \sigma_{\alpha}(n) = x^{\alpha}\exp\left(\left(\frac{1}{1-\alpha}+o(1)\right)\frac{(\ln x)^{1-\alpha}}{\ln \ln x}\right)
    \end{equation}
    hence
    \begin{equation}
        E_{\alpha}(x) = \Omega\left(x^{\alpha}\exp\left(\left(\frac{1}{1-\alpha}+o(1)\right)\frac{(\ln x)^{1-\alpha}}{\ln \ln x}\right)\right)
    \end{equation}
\end{theorem}
and here we give a short proof.

\begin{proof}
We have
    \begin{equation}
        \begin{aligned}
            &\sigma_{\alpha}(n) = n^{\alpha} \sigma_{-\alpha}(n) = n^{\alpha} \prod_{p|n} (1+p^{-\alpha}+O(p^{-2\alpha})) \\
            &= n^{\alpha} \exp\left(\sum\limits_{p|n} p^{-\alpha}+O\left(\sum\limits_{p|n} p^{-2\alpha}\right)\right)
        \end{aligned}
    \end{equation}

Note that 
\begin{equation}
    \sup_{n\leq x} \sum_{p|n} p^{-\alpha} = \sum_{p|P(Y_x)} p^{-\alpha}
\end{equation}

where $P(y) = \prod\limits_{p\leq y} p$ and $Y_x$ is the largest positive integer that satisfies $P(Y_x) \leq x$. 

Hence 
\begin{equation}
    \sup_{n\leq x} \sigma_{\alpha}(n) =x^{\alpha}\exp\left(\sum_{p|P(Y_x)} p^{-\alpha}+O\left(\sum_{p|P(Y_x)} p^{-2\alpha}\right)\right)
\end{equation}

If $\beta\geq1$, $\sum\limits_{p|P(Y_x)} p^{-\beta}$ is bounded.

If $0<\beta < 1$, by Prime Number Theorem and Abel summation, we have
\begin{equation}
    \sum\limits_{p|P(Y_x)} p^{-\beta} \sim \frac{Y_x^{1-\beta}}{(1-\beta) \ln Y_x} \sim \frac{(\ln x)^{1-\alpha}}{(1-\beta) \ln \ln x}
\end{equation}

Applying the above observations for $\beta=\alpha$ and $2\alpha$, we obtain
\begin{equation}
    \sup_{n\leq x} \sigma_{\alpha}(n) = x^{\alpha} \exp\left(\frac{(\ln x)^{1-\alpha}}{ \ln \ln x} \left(\frac{1}{1-\alpha}+o(1)\right)\right)
\end{equation}
as desired.
\end{proof}

\section{Proof of Theorems}

Consider $\sum\limits_{n=1}^{\infty} \frac{\sigma_{\alpha}(n)}{n^s} = \zeta(s) \zeta(s-\alpha)$.  By Mellin's Inversion Formula, we have 
\begin{equation} \label{exp_sigma_alpha}
    \sum_{n \leq x} \sigma_{\alpha}(n) = \frac{1}{2\pi i} \lim_{T\rightarrow +\infty} \int_{c-iT}^{c+iT}f(s)\frac{x^s}{s}ds
\end{equation}
where $f(s) = \sum\limits_{n=1}^{\infty} \frac{\sigma_{\alpha}(n)}{n^s}=\zeta(s)\zeta(s-\alpha)$ and we choose $c=1+\alpha+\varepsilon$ with $\varepsilon>0$.

Nevertheless, we wish to remove the infinity in the expression and approximate $\sum\limits_{n\leq x} \sigma_{\alpha}(n)$ with $I_1(x) = \frac{1}{2\pi i} \int\limits_{c-iT}^{c+iT}f(s)\frac{x^s}{s}ds$.

Take $T\geq 1$, then
\begin{equation} \label{3a}
    \begin{aligned} 
        I_1(x)&= \frac{1}{2\pi i} \int_{c-iT}^{c+iT}f(s)\frac{x^s}{s}ds \\ 
        &= \frac{1}{2\pi i} \int_{c-iT}^{c+iT}\left(\sum_{n=1}^{\infty} \frac{\sigma_{\alpha}(n)}{n^s}\right)\frac{x^s}{s}ds \\
        &= \sum_{n=1}^{\infty} \sigma_{\alpha}(n) \frac{1}{2\pi i}\int_{c-iT}^{c+iT} \frac{x^s}{sn^s}ds \\
    \end{aligned}
\end{equation}

Notice that now we can apply Perron's formula \cite{tenenbaum2015introduction} to (\ref{3a}):
\begin{equation} \label{perron}
    \frac{1}{2\pi i}\int_{c-iT}^{c+iT} \frac{a^s}{s} ds = 
        \begin{cases}
         1+O(\frac{a^c}{T\ln a})& \text{ if } a>1 \\
         O(\frac{a^c}{T|\ln a|})& \text{ if } a<1 
        \end{cases}
\end{equation}

and we have 
\begin{equation} \label{3b}
    \begin{aligned} 
        I_1(x)&= \sum_{n\geq 1} \sigma_{\alpha}(n) \mathbbm{1}_{n<x} +O\left(\frac{x^c}{Tn^c |\ln{\frac{x}{n}}|}\right) \\
        &= \sum_{n\leq x} \sigma_{\alpha}(n) + O\left(x^c \sum_{n\geq 1} \frac{\sigma_{\alpha}(n)}{Tn^c|\ln{\frac{x}{n}}|}\right)
    \end{aligned}
\end{equation}

Theorem 3 implies the following bound for $\sigma_\alpha(n)$:

\begin{lemma} \label{sigma_alpha_bound}
    \begin{equation}
        \sigma_{\alpha}(x) = O(x^{\alpha+o(1)})
    \end{equation}
\end{lemma}

We have that if $n \leq \frac{x}{2}$ or $n \geq 2x$, then $|\ln{\frac{x}{n}}| \geq \ln 2$; otherwise, $\ln{\frac{x}{n}} \geq \frac{|n-x|}{x}$. Plugging in these bounds into (\ref{3b}) and utilizing and Lemma (\ref{sigma_alpha_bound}), we get
\begin{equation} \label{3c}
    \begin{aligned} 
        O\left(x^c\sum_{n\geq 1} \frac{\sigma_{\alpha}(n)}{Tn^c|\ln{\frac{x}{n}}|}\right) &= O\left(x^c\sum_{n\leq \frac{x}{2} \text{ or } n\geq 2x} \frac{\sigma_{\alpha}(n)}{Tn^c|\ln{\frac{x}{n}}|}\right) + O\left(x^c\sum_{\frac{x}{2} < n < 2x} \frac{\sigma_{\alpha}(n)}{Tn^c|\ln{\frac{x}{n}}|}\right) \\
        &= O\left(x^c\sum_{n \geq 1} \frac{n^{\alpha+o(1)}}{Tn^c}\right) + O\left(x^{c+1}\sum_{\frac{x}{2} < n < 2x} \frac{n^{\alpha+o(1)}}{Tn^c |n-x|}\right) \\
        &=O\left(\frac{x^c}{T}\sum_{n\geq 1}n^{\alpha-c+o(1)}\right) + O\left(\frac{x^{\alpha+1+o(1)}}{T}\right)
    \end{aligned}
\end{equation}

For the first term, further notice that
\begin{equation} \label{3d}
    \begin{aligned} 
       \sum_{n\geq 1}n^{\alpha-c+o(1)} < 1+\int_{1}^{\infty} n^{\alpha-c+o(1)} dn = 1-\frac{1}{1+\alpha-c+o(1)} < 1+\frac{2}{\varepsilon}
    \end{aligned}
\end{equation}

So we have 

\begin{equation} \label{3e}
     O\left(x^c\sum_{n\geq 1} \frac{\sigma_{\alpha}(n)}{Tn^c|\ln{\frac{x}{n}}|}\right) = O\left(\frac{x^{c+o(1)}}{T}\right) + O\left(\frac{x^{\alpha+1+o(1)}}{T}\right) = O\left(\frac{x^{c+o(1)}}{T}\right)
\end{equation}

Plugging back to (\ref{3b}), we get our desired estimation:

\begin{equation} \label{sum_finite_est}
    \sum_{n\leq x} \sigma_{\alpha}(n) = I_1(x)+ O\left(\frac{x^{c+o(1)}}{T}\right)
\end{equation}

Our following work is focused on determining the bound of error term in $I_1(x)$. Now that the integral involved in $I_1(x)$ is along the real axis $\sigma = c = 1+\alpha+\varepsilon$, we can shift the contour to $\sigma = 1+\alpha-c = -\varepsilon$ through Cauchy's residue theorem and then apply the functional equation of Riemann zeta function.
\begin{equation} \label{3f}
    I_1(x) =-I_2(x)-I_3(x)-I_4(x)+ \nu(x)
\end{equation}
where
\begin{equation} \label{I_2}
    I_2(x) = \frac{1}{2\pi i} \int_{c+iT}^{1+\alpha-c+iT} f(s)\frac{x^s}{s}ds
\end{equation}

\begin{equation} \label{I_3}
    I_3(x) = \frac{1}{2\pi i} \int_{1+\alpha-c+iT}^{1+\alpha-c-iT} f(s)\frac{x^s}{s}ds
\end{equation}

\begin{equation} \label{I_4}
    I_4(x) = \frac{1}{2\pi i} \int_{1+\alpha-c-iT}^{c-iT} f(s)\frac{x^s}{s}ds
\end{equation}
and
\begin{equation} \label{residue}
    \begin{aligned}
        \nu(x) &= \mathrm{Res}\left(f(s)\frac{x^s}{s},\alpha+1\right) + \mathrm{Res}\left(f(s)\frac{x^s}{s},0\right) + \mathrm{Res}\left(f(s)\frac{x^s}{s},1)\right) \\
        &= \frac{x^{1+\alpha}}{1+\alpha}\zeta(1+\alpha)-\frac{\zeta(-\alpha)}{2}+x\zeta(1-\alpha)
    \end{aligned}
\end{equation}

We first estimate $I_2$ and $I_4$.

\begin{equation} \label{I_2_est}
    \begin{aligned}
        I_2(x) &= O\left( \int_{1+\alpha-c+iT}^{c+iT} f(s)\frac{x^s}{s}ds\right)\\
        &= O\left(\int_{1+\alpha-c}^{c}|\zeta(\sigma+iT)\zeta(\sigma-\alpha+iT)|\frac{x^{\sigma}}{T}d\sigma\right)\\
        &= O\left(\int_{1+\alpha-c}^{c} T^{\mu(\sigma)+\mu(\sigma-\alpha)+\varepsilon-1}x^{\sigma}d\sigma\right)
    \end{aligned}
\end{equation}

where the last step comes from the following bound for $|\zeta(\sigma+it)|$ from \cite{titchmarsh1986theory}:

\begin{lemma} \label{zeta_est}
    \begin{equation} 
    |\zeta(\sigma+it)| \ll t^{\mu(\sigma)} \ln t,
    \mu(\sigma) = 
        \begin{cases}
         0& \text{ if } \sigma\geq1 \\
         \frac{1-\sigma}{3}& \text{ if } \frac{1}{2}\leq\sigma<1 \\
         \frac{1}{2}-\frac{2\sigma}{3}& \text{ if } \sigma<\frac{1}{2}
        \end{cases}
\end{equation}
\end{lemma}

After integrating $T^{\mu(\sigma)+\mu(\sigma-\alpha)+\varepsilon-1}x^{\sigma}$ by partitioning $[1+\alpha-c, c]$ into separate intervals corresponding to different values of $\mu(\sigma)$ and $\mu(\sigma-\alpha)$, we get

\begin{equation} \label{I_2_est_final}
    I_2(x) = O\left(\frac{x^{1+\alpha+\varepsilon}}{T^{1-\varepsilon}}\right)
\end{equation}

Similar results can be deduced for $I_4(x)$:
\begin{equation} \label{I_4_est_final}
    I_4(x) = O\left(\frac{x^{1+\alpha+\varepsilon}}{T^{1-\varepsilon}}\right)
\end{equation}

Now our task is to estimate $I_3$. We can apply the functional equation of zeta function on $f(s)$ in the integral. 

\begin{lemma} \label{func_eq_zeta}
    (Functional equation of Riemann zeta function) 
    \begin{equation}
        \zeta(s)\pi^{-\frac{s}{2}} \Gamma(\frac{s}{2}) = \zeta(1-s)\pi^{-\frac{1-s}{2}}\Gamma(\frac{1-s}{2})
    \end{equation}
\end{lemma}

\begin{corollary} \label{func_eq_prod}
    \begin{equation}
        \zeta(s)\zeta(s-\alpha)\Gamma(\frac{s}{2})\Gamma(\frac{s-\alpha}{2}) = \zeta(1-s)\zeta(1+\alpha-s)\pi^{2s-\alpha-1}\Gamma(\frac{1-s}{2})\Gamma(\frac{1+\alpha-s}{2})
    \end{equation}
\end{corollary}

We define $\gamma(s)$ to be the factor in the functional equation of $f(s)$ by applying Corollary \ref{func_eq_prod}:

\begin{equation} \label{func_eq_f}
    f(1+\alpha-s) = \gamma(s)f(s), \gamma(s)=\pi^{1+\alpha-2s}\frac{\Gamma(\frac{s}{2})\Gamma(\frac{s-\alpha}{2})}{\Gamma(\frac{1-s}{2})\Gamma(\frac{1+\alpha-s}{2})}
\end{equation}

By utilizing (\ref{func_eq_f}) on $I_3(x)$, we get
\begin{equation}
    \begin{aligned}
        I_3(x) &= \frac{1}{2\pi i}\int_{1+\alpha-c+iT}^{1+\alpha-c-iT} f(s)\frac{x^s}{s} ds \\
        &= \frac{1}{2\pi i} x^{1+\alpha} \int_{c+iT}^{c-iT} f(s)\gamma(s)\frac{x^{-s}}{1+\alpha-s}ds \\
        &= x^{1+\alpha}\sum_{m\geq 1} \sigma_{\alpha}(m)c_T(mx)
    \end{aligned}
\end{equation}

where $c_T(y) = \frac{1}{2\pi i}\int\limits_{c+iT}^{c-iT} \frac{y^{-s}\gamma(s)}{1+\alpha-s}ds$.

Now the most difficult task is to estimate $c_T(y)$.

To estimate $\gamma(s)$, we use the following standard bound for $\Gamma(\sigma+it)$ from Stirling's formula when $t$ is large:

\begin{equation} \label{Gamma_est}
    \Gamma(\sigma+it)=\sqrt{2\pi}(it)^{\sigma-\frac{1}{2}} e^{-\frac{\pi}{2}t} \left(\frac{t}{e}\right)^{it}(1+O(t^{-1}))
\end{equation}

By plugging in this into $\gamma(s)$, we have
\begin{equation} \label{self_gamma_est}
    \gamma(\sigma+it) = \left(\frac{t}{2\pi}\right)^{2\sigma-\alpha-1} e^{4\pi i \theta(t)} \left(1+O\left(\frac{1}{t}\right)\right), \theta(t) = \frac{t}{2\pi} \ln\left(\frac{t}{2\pi}\right) -\frac{t}{2\pi} -\frac{1}{8}
\end{equation}

Then
\begin{equation}
    \begin{aligned}
        |c_T(y)| &= \mathfrak{Re} \left(\frac{1}{\pi} \int_{c}^{c+iT} y^{-s}\frac{\gamma(s)}{1+\alpha-s}ds\right) \\ 
        &= \mathfrak{Re}\left( \frac{1}{\pi}\left(\int_{c}^{c+T} y^{-s}\frac{\gamma(s)}{1+\alpha-s}ds + \int_{c+T}^{c+iT} y^{-s}\frac{\gamma(s)}{1+\alpha-s}ds\right)\right) \\
        &= O(y^{-c}) + \mathfrak{Re}\int_{1}^{T} e^{4\pi i\theta(t)-it\ln y} \left(\frac{t}{2\pi}\right)^{2c-\alpha-1} \left(-\frac{1}{t}\right)(1+O(t^{-1}))dt
    \end{aligned}
\end{equation}

From \cite{ivic2012riemann}, the following inequality holds:

\begin{equation}
    \int_{J} e^{2\pi i g(t)}h(t) dt \ll \frac{\max{|h(t)|}}{\delta} \text{ if } |g'(t)| > \delta \text{ for all } t\in J
\end{equation}

In our case $g(t) = 2\theta(t)-\frac{t\ln y}{2\pi}$, $g'(t) = \frac{1}{\pi} \ln(\frac{t}{2\pi}) - \frac{\ln y}{2\pi}$. When $y > 2T^2$, we have that $g'(t)< U(T)$ where $U(T)$ is a function of $T$.

Then $|c_T(y)| = O(y^{-c}) + O(y^{-c}T^{2c-\alpha-2}) = O(y^{-c}T^{2c-\alpha-2})$ when $y> 2T^2$.

When $y \leq 2T^2$, we move the contour of integration to $b=1+\frac{\alpha}{2}$.

Since there are no singularities in the contour of integration, we have

\begin{equation}
    \begin{aligned}
        c_T(y) &= \frac{1}{2\pi i} \left(\int_{b-iT}^{b+iT} \frac{\gamma(s)}{\alpha+1-s}y^{-s}ds + \int_{b+iT}^{c+iT} \frac{\gamma(s)}{\alpha+1-s}y^{-s}ds + \int_{c-iT}^{b-iT} \frac{\gamma(s)}{\alpha+1-s}y^{-s}ds\right) \\
        &= \frac{1}{2\pi i} \int_{b-iT}^{b+iT} \frac{\gamma(s)}{\alpha+1-s}y^{-s}ds + O\left(\int_{b}^{c} T^{2\sigma-\alpha-2}y^{-\sigma}d{\sigma}\right) \\
        &= \frac{1}{2\pi i} \int_{b-iT}^{b+iT} \frac{\gamma(s)}{\alpha+1-s}y^{-s}ds + O\left(T^{-\alpha-2} \left(\left(\frac{T^2}{y}\right)^c + \left(\frac{T^2}{y}\right)^b\right)\right) \\
        &= \frac{1}{\pi} \mathfrak{Re}\int_{1}^{T} \frac{\gamma(\sigma+it)}{\alpha+1-b-it}y^{-b-it}dt + O(T^{\alpha+2\varepsilon}y^{-c}+y^{-b})
    \end{aligned}
\end{equation}

Using (\ref{self_gamma_est}) and $\frac{1}{a+1-b-it} = \frac{a+1-b+it}{(a+1-b)^2+t^2} = O\left(\frac{1}{t}\right)$, 

\begin{equation}
    \begin{aligned}
        c_T(y) &= \frac{i}{2\pi^2} y^{-b}\mathfrak{Re}\left(\int_{1}^{T}e^{4\pi i \theta(t)}y^{-it} O((1+t^{-1})^2)dt\right) + O(T^{\alpha+2\varepsilon}y^{-c}+y^{-b}) \\
        &= \frac{i}{2\pi^2} y^{-b}\mathfrak{Re}\left(\int_{1}^{T}e^{4\pi i \theta(t)}y^{-it} dt\right) +  O(T^{\alpha+2\varepsilon}y^{-c}+y^{-b}\ln T)
    \end{aligned}
\end{equation}

From Lemma 2.1 in \cite{friedlander2005summation}:

\begin{lemma}
    If $2 \leq z \leq 2T$, then
    \begin{equation}
        \int_{1}^T \left(\frac{t}{ez}\right)^{imt} dt = \sqrt{\frac{2\pi z}{m}} e^{\frac{\pi i}{4} -imz} + O\left(\left(\frac{1}{\sqrt{T}}+\left| \log \frac{z}{T}\right|^{-1}\right)\right)
    \end{equation}
\end{lemma}

We have
\begin{equation}
    \begin{aligned}
        \int_1^{T} e^{4\pi i \theta(t)} y^{-it} dt &= -i \int_1^{T} \left(\frac{t}{2\pi e \sqrt{y}}\right)^{2it}dt \\
        &= 
        \begin{cases}
            -\sqrt{2}\pi i y^{\frac{1}{4}} e^{\frac{\pi i}{4}-4\pi \sqrt{y}i} + O\left(\left(\frac{1}{\sqrt{T}}+\ln \frac{T}{2\pi \sqrt{y}}\right)^{-1}\right) & y \leq \frac{T^2}{4\pi^2}\\
            O\left(\left(\frac{1}{\sqrt{T}}+\ln \frac{2\pi \sqrt{y}}{T}\right)^{-1}\right) & \frac{T^2}{4\pi^2} < y \leq 2T^2
        \end{cases}
    \end{aligned}
\end{equation}

Hence, 
\begin{equation}
    \begin{aligned}
        c_T(y) = 
        \begin{cases} \frac{1}{\sqrt{2}\pi}y^{-\frac{3}{4}-\frac{\alpha}{2}}\cos(4\pi\sqrt{y}-\frac{\pi}{4}) + O\left(y^{-1-\frac{\alpha}{2}}\left(\frac{1}{\sqrt{T}}+\ln \frac{T}{2\pi \sqrt{y}}\right)^{-1}\right) \\ + O(T^{\alpha+2\varepsilon}y^{-c}+y^{-1-\frac{\alpha}{2}}\ln T) & y\leq \frac{T^2}{4\pi^2} \\
        O\left(y^{-1-\frac{\alpha}{2}}\left(\frac{1}{\sqrt{T}}+\ln \frac{2\pi \sqrt{y}}{T}\right)^{-1}\right) + O(T^{\alpha+2\varepsilon}y^{-c}+y^{-1-\frac{\alpha}{2}}\ln T) & \frac{T^2}{4\pi^2} < y \leq 2T^2\\
        O(y^{-c}T^{2c-\alpha-2}) & y> 2T^2
        \end{cases}
    \end{aligned}
\end{equation}
Let $y = mx$.
\begin{equation} \label{I_3_est_final}
    \begin{aligned}
        I_3(x) &= x^{1+\alpha}\sum_{m\geq 1} \sigma_{\alpha}(m)c_T(mx) \\
            &= x^{1+\alpha}\left(\sum_{mx \leq \frac{T^2}{4\pi^2}} \sigma_{\alpha}(m)\left(\frac{1}{\sqrt{2}\pi}y^{-\frac{3}{4}-\frac{\alpha}{2}}\cos\left(4\pi\sqrt{y}-\frac{\pi}{4}\right)\right)\right. \\
            &+ O\left(\left(\frac{T^2}{x}\right)^{\frac{\alpha}{2}+o(1)}x^{-1-\frac{\alpha}{2}}+\left(\frac{T^2}{x}\right)^{\frac{\alpha}{2}+o(1)-1}\sqrt{T}x^{-1-\frac{\alpha}{2}}\right) \\ & \left . +O\left(\left(\frac{T^2}{x}\right)^{\alpha+1-c+o(1)}T^{\alpha+2\varepsilon}x^{-c}\right)+O\left(\left(\frac{T^2}{x}\right)^{\frac{\alpha}{2}+o(1)}x^{-1-\frac{\alpha}{2}}\ln T\right) + O(T^{2c-\alpha-2}x^{-c})\vphantom{\sum_{mx \leq \frac{T^2}{4\pi^2}} \sigma_{\alpha}(m)(\frac{1}{\sqrt{2}\pi}y^{-\frac{3}{4}-\frac{\alpha}{2}}e^{\frac{\pi i}{4}-4\pi \sqrt{y}i})}\right) \\
            &= x^{1+\alpha}\sum_{mx \leq \frac{T^2}{4\pi^2}} \sigma_{\alpha}(m)\left(\frac{1}{\sqrt{2}\pi}y^{-\frac{3}{4}-\frac{\alpha}{2}}\cos\left(4\pi\sqrt{y}-\frac{\pi}{4}\right)\right) + O(T^{\alpha-\frac{3}{2}+o(1)}x^{1-o(1)})\\ &+ O(T^{\alpha+o(1)}x^{-o(1)}) + O(T^{2c-\alpha-2}x^{1+\alpha-c})
            \\ 
            &= x^{1+\alpha}\sum_{mx \leq \frac{T^2}{4\pi^2}} \sigma_{\alpha}(m)\left(\frac{1}{\sqrt{2}\pi}y^{-\frac{3}{4}-\frac{\alpha}{2}}\cos\left(4\pi\sqrt{y}-\frac{\pi}{4}\right)\right) + O(T^{\alpha-\frac{3}{2}+o(1)}x^{1-o(1)})\\ &+ O(T^{\alpha+o(1)} x^{-o(1)}) + O(T^{\alpha+2\varepsilon}x^{-\varepsilon})
    \end{aligned}
\end{equation}

Summing up (\ref{sum_finite_est}) (\ref{3f}) (\ref{residue}) (\ref{I_2_est_final}) (\ref{I_4_est_final}) (\ref{I_3_est_final}), we have
\begin{equation}
    \begin{aligned}
        \sum_{n\leq x} \sigma_{\alpha}(n) &= \frac{x^{1+\alpha}}{1+\alpha}\zeta(1+\alpha)-\frac{\zeta(-\alpha)}{2}+x\zeta(1-\alpha) \\&+ x^{1+\alpha}\sum_{mx \leq \frac{T^2}{4\pi^2}} \sigma_{\alpha}(m)\left(\frac{1}{\sqrt{2}\pi}y^{-\frac{3}{4}-\frac{\alpha}{2}}\cos\left(4\pi\sqrt{y}-\frac{\pi}{4}\right)\right) \\ &+ O(T^{\varepsilon-1}x^{1+\alpha+\varepsilon}) + O(T^{\alpha-\frac{3}{2}+o(1)}x^{1-o(1)}) + O(T^{\alpha+o(1)} x^{-o(1)}) + O(T^{\alpha+2\varepsilon}x^{-\varepsilon}) \\
        &= \frac{x^{1+\alpha}}{1+\alpha}\zeta(1+\alpha)-\frac{\zeta(-\alpha)}{2}+x\zeta(1-\alpha) \\&+ x^{\frac{1}{4}+\frac{\alpha}{2}}\sum_{mx \leq \frac{T^2}{4\pi^2}} \sigma_{\alpha}(m)\left(\frac{1}{\sqrt{2}\pi}m^{-\frac{3}{4}-\frac{\alpha}{2}}\cos\left(4\pi\sqrt{y}-\frac{\pi}{4}\right)\right) \\
        &+ O\left(x^{o(1)}\left(\frac{x^{1+\alpha}}{T}+T^{\alpha}\right)\right)
    \end{aligned}
\end{equation}

This gives

\begin{equation}
    \begin{aligned}
        E_{\alpha}(x)&=O\left(x^{\frac{1}{4}+\frac{\alpha}{2}}\left(\frac{T^2}{x}\right)^{\frac{1}{4}+\frac{\alpha}{2}}\right)+O\left(x^{o(1)}\left(\frac{x^{1+\alpha}}{T}+T^{\alpha}\right)\right)\\
        &= O(T^{\frac{1}{2}+\alpha})+O\left(\frac{x^{1+\alpha+o(1)}}{T}\right)
    \end{aligned}
\end{equation}
Pick $T=x^{\frac{2+2\alpha}{3+2\alpha}}$. Then $E_\alpha(x)=O\left(x^{\frac{2\alpha^2+3\alpha+1}{2\alpha+3}+\varepsilon}\right)$, proving Corollary \ref{intermediate_bound}.

Now fix $x^{\frac{3}{4}+\frac{\alpha}{2}+\varepsilon} \ll T \ll x^{\frac{1}{4\alpha}+\frac{1}{2}-\varepsilon}$. Then $O(\frac{x^{1+\alpha+O(1)}}{T}) = O(\frac{x^{1+\alpha}}{x^{\frac{3}{4}+\frac{\alpha}{2}}}) = O(x^{\frac{1}{4}+\frac{\alpha}{2}})$, $O(x^{o(1)}T^{\alpha}) = O((x^{\frac{1}{4\alpha}+\frac{1}{2}})^{\alpha}) = O(x^{\frac{1}{4}+\frac{\alpha}{2}})$. So 

\begin{equation}
    E_{\alpha}(x) = x^{\frac{1}{4}+\frac{\alpha}{2}}\sum_{mx \leq \frac{T^2}{4\pi^2}} \sigma_{\alpha}(m)\left(\frac{1}{\sqrt{2}\pi}m^{-\frac{3}{4}-\frac{\alpha}{2}}\cos\left(4\pi\sqrt{y}-\frac{\pi}{4}\right)\right) + O(x^{\frac{1}{4}+\frac{\alpha}{2}})
\end{equation}

In addition, let  $x^{\frac{1}{2}+\alpha+\varepsilon} \ll N \ll x^{\frac{1}{2\alpha}-\varepsilon}$. Define $$F_{\alpha}(x, N) = \sum_{n\leq N} \frac{\sigma_{\alpha}(n)}{n^{\frac{3}{4}+\frac{\alpha}{2}}} \cos(4\pi \sqrt{nx}-\frac{\pi}{4})$$ under this constraint on $N$. Our goal is to find $x \in (N^{\frac{2}{1+2\alpha}-2\delta}, N^{\frac{2}{1+2\alpha}-\delta})$ ($\delta>0$) such that $F_\alpha(x, N) = \Omega((\ln x)^{\frac{1}{4}+\frac{\alpha}{2}})$.

Now we apply Soundararajan's Lemma \cite{soundararajan2003omega} on $F_{\alpha}(x^2,N)$ to give an $\Omega$-theorem. The lemma is summarized in the following statement:

\begin{lemma} \label{lem:Sound's}
    (Soundararajan's Lemma) Let $(f(n))_{n=1}^{\infty}$ and $(\lambda _n)_{n=1}^{\infty}$ be sequences of nonnegative real numbers, $(\lambda_n)$ non-decreasing, and $\sum\limits_{n=1}^{\infty}f(n) < \infty$. Let $L \geq 2$ and $Y$ be integers, $\gamma \in \mathbb{R}$ fixed. Let 
    \begin{equation}
        S(t) = \sum_{n=1}^{\infty}f(n) \cos(2\pi \lambda_n t + \gamma)
    \end{equation}
    Suppose $\mathcal{M}$ is a finite set of positive integers such that $\{\lambda_m:m\in M\} \subset [\frac{\lambda_Y}{2}, \frac{3\lambda_Y}{2}]$. Then for any real $X \geq 2$, there exists some $x\in [\frac{X}{2}, (6L)^{\mid \mathcal{M} \mid +1}X]$ such that 
    \begin{equation}
        \mid S(x) \mid \geq \frac{1}{8}\sum_{m\in \mathcal{M}}f(m)-\frac{1}{L-1}\sum_{\{n \mid \lambda_n \leq 2\lambda_Y \}}f(n) -\frac{4}{\pi^2 X \lambda_Y} \sum_{n=1}^{\infty} f(n)
    \end{equation}
\end{lemma}

In our application, $S(x) = F_{\alpha} (x^2, N)$. We have $F_{\alpha}(x^2, N) = 
\sum_{n=1}^{\infty}f(n) \cos(2 \pi \lambda_n x + \beta)$, where $\lambda_n = 2\sqrt{n}$, $\beta = -\frac{\pi}{4}$,

\begin{equation}
    f(n) = \left\{ \begin{array}{cl}
\frac{\sigma_{\alpha}(n)}{n^{\frac{\alpha}{2}+\frac{1}{4}}} &  \ n \leq N \\
0 &  \ \text{otherwise}
\end{array} \right.
\end{equation}.

Take $\mathcal{M} = [\frac{T}{4}, \frac{3T}{4}]$.

We require $\frac{X}{2} > N^{\frac{1}{1+2\alpha}-\delta}$ and $(6L)^{\mid \mathcal{M} \mid +1}X < N^{\frac{1}{1+2\alpha}-\frac{\delta}{2}}$. This gives $(6L)^{\mid \mathcal{M} \mid +1} = (6L)^{2T+1} < \frac{N^{\frac{\delta}{2}}}{2}$. We take $L$ to be a large constant (say $L=10^8$). Then $T < \frac{\delta \ln N}{4  \ln(6 \cdot 10^8)}$.

Take $T = \frac{\delta \ln N}{8  \ln(6 \cdot 10^8)}$.

\begin{equation} \label{part_1}
    \sum_{m\in \mathcal{M}} f(m) \geq \sum_{\frac{T}{4} \leq m \leq \frac{3T}{4}} m^{-\frac{3}{4}+\frac{\alpha}{2}} \gg T^{\frac{1}{4}+\frac{\alpha}{2}}
\end{equation}

\begin{equation} \label{part_2}
    \sum_{\{n \mid \lambda_n \leq 2\lambda_Y \}}f(n) \ll T^{\frac{1}{4}+\frac{\alpha}{2}}
\end{equation}

\begin{equation}
    \frac{1}{X\lambda_T} \ll \frac{1}{X\sqrt{T}} \ll \frac{1}{\sqrt{\ln N}} N^{\delta-\frac{1}{1+2\alpha}}
\end{equation}

\begin{equation}
    \sum_{n=1}^{\infty} f(n) = \sum_{n \leq N} \frac{\sigma_{\alpha}(n)}{n^{\frac{3}{4}+\frac{\alpha}{2}}} \ll N^{\frac{1}{4}+\frac{\alpha}{2}}
\end{equation}

So 
\begin{equation} \label{part_3}
    \frac{1}{X\lambda_T}\sum_{n=1}^{\infty} f(n) \ll \frac{N^{\delta}}{\sqrt{\ln N}} N^{\frac{1}{4}+\frac{\alpha}{2}-\frac{1}{1+2\alpha}}.
\end{equation}

Note that $N \ll x^{\frac{1}{2}}$. In order to confine $\frac{1}{X\lambda_T}\sum\limits_{n=1}^{\infty} f(n)$ at the order of $(\ln x)^{\frac{1}{4}+\frac{\alpha}{2}}$, we need that $\frac{1}{4} + \frac{\alpha}{2}-\frac{1}{1+2\alpha} < 0$. Solving the inequality gives $\alpha <\frac{1}{2}$, which perfectly satisfies the assumption in Theorem \ref{thm:main}.

Plugging (\ref{part_1}) (\ref{part_2}) (\ref{part_3}) into Lemma \ref{lem:Sound's}, we get $\mid F_{\alpha}(x^2, N) \mid \geq T^{\frac{1}{4}+\frac{\alpha}{2}} = \Omega( (\ln N)^{\frac{1}{4}+\frac{\alpha}{2}}) = \Omega((\ln x^2)^{\frac{1}{4}+\frac{\alpha}{2}})) $, thus fulfilling our goal.

So $E_{\alpha}(x) = \Omega((x\ln x)^{\frac{1}{4}+\frac{\alpha}{2}})$.

\section{Concluding remarks}

It is highly probable that our bound $\Omega((x\ln x)^{1/4+\alpha/2})$ is not optimal. However, trying to improve it using Soundararajan's method meets substantial obstacles, mainly because the function
$$
\omega_\alpha(n)=\sum_{p\mid n}p^{-\alpha}
$$
behaves very differently from $\omega=\omega_0$, which is the number of distinct prime divisors. To be more precise, $\omega_\alpha$ is concentrated around its average, which is $O(1)$ in this case, hence the sets $\mathcal M\subset [1,N]$ such that $\sigma_\alpha(m)$ is large for all $m\in \mathcal M$ are much smaller than the analogous sets for $d=\sigma_0$.

Results like Theorem 2 are a consequence of a summation formula analogous to the Voronoi summation formula. The corresponding kernel of an integral transform takes form
$$
W_{\alpha}(y) = \lim\limits_{T\rightarrow \infty} \frac{1}{2\pi i} \int\limits_{2-iT}^{2+iT} y^{-s} \frac{\Gamma(\frac{s}{2}) \Gamma(\frac{s-\alpha}{2})}{(\alpha+1-s) \Gamma(\frac{1-s}{2})\Gamma(\frac{1+\alpha-s}{2})}  ds.
$$
One can prove that similarly to the original Voronoi summation, the kernel can be expressed in terms of Bessel functions, namely
$$
 W_{\alpha}(y) = \frac{1}{2\sin(\frac{\pi \alpha}{2}) y^{\frac{1+\alpha}{2}}}(I_{-1-\alpha}(4\sqrt{y})-I_{\alpha+1}(4\sqrt{y})-J_{-\alpha-1}(4\sqrt{y})-J_{\alpha+1}(4\sqrt{y})).
$$
This function interpolates between the kernels for $d(n)$ and $\sigma(n)$ for $\alpha \in [0,1]$.

\bibliographystyle{alpha}
\bibliography{refs}

\end{document}